\documentclass[11pt]{amsart}

\usepackage{amsfonts, amstext, amsmath, amsthm, amscd, amssymb, upgreek}
\usepackage{graphicx, color,  subfigure, wrapfig, overpic, url, wasysym}
\usepackage[dvipsnames]{xcolor}

\newtheorem{question}{Question}

\newcommand{\asplice}{\ensuremath{\makebox[0.3cm][c]{\raisebox{-0.3ex}{\rotatebox{90}{$\asymp$}}}}}

\begin{document}

\title{Turaev Surfaces} 

\author{Seungwon Kim}
\address{National Institute for Mathematical Sciences, Daejeon, South Korea}

\author{Ilya Kofman}
\address{Department of Mathematics, College of Staten Island \& The Graduate Center, City University of New York, New York, USA}


\thanks{An expository article written for \emph{A Concise Encyclopedia of Knot Theory}, to appear.}

\maketitle

\section{Introduction}\label{sec:intro}

The two most famous knot invariants, the Alexander polynomial (1923)
and the Jones polynomial (1984), mark paradigm shifts in knot theory.
After each polynomial was discovered, associating new structures to
knot diagrams played a key role in understanding its basic properties.
Using the Seifert surface, Seifert showed, for example, how to obtain
knots with a given Alexander polynomial.  For the Jones polynomial,
even the simplest version of that problem remains open: Does there
exist a non-trivial knot with trivial Jones polynomial?

Kauffman gave a state sum for the Jones polynomial, with terms for
each of the $2^c$ states of a link diagram with $c$ crossings.  Turaev
constructed a closed orientable surface from any pair of dual states
with opposite markers.

For the Jones polynomial, the Turaev surface is a rough analog to the
Seifert surface for the Alexander polynomial.  For a given knot
diagram, the Seifert genus and the Turaev genus are computed by
separate algorithms to obtain each surface from the diagram.  The
invariants for a given knot $K$ are defined as the minimum genera
among all the respective surfaces for $K$.  The Seifert genus is a
topological measure of how far a given knot is from being unknotted.
The Turaev genus is a topological measure of how far a given knot is
from being alternating.  (See \cite{Lowrance-knot_encyclopedia}, which
discusses alternating distances.)  For any alternating diagram,
Seifert's algorithm produces the minimal genus Seifert surface.  For
any adequate diagram, Turaev's algorithm produces the minimal genus
Turaev surface.  Extending the analogy, we can determine the Alexander
polynomial and the Jones polynomial of $K$ from associated algebraic
structures on the respective surfaces of $K$: the Seifert matrix for
the Alexander polynomial, and the $A$--ribbon graph on the Turaev
surface for the Jones polynomial.

The analogy is historical, as well.  Like the Seifert surface for the
Alexander polynomial, the Turaev surface was constructed to prove a
fundamental conjecture related to the Jones polynomial.  In the
1880's, Tait conjectured that an alternating link always has an
alternating diagram that has minimal crossing number among all
diagrams for that link.  A proof had to wait about a century until the
Jones polynomial led to several new ideas used to prove Tait's
Conjecture \cite{kauffman, murasugi, Thistlethwaite}. Turaev's later
proof in \cite{turaev} introduced Turaev surfaces and prompted
interest in studying their properties.

\section{What is the Turaev surface?}\label{sec:def}

Let $D$ be the diagram of a link $L$ drawn on $S^2$.  For any crossing
\includegraphics[height=0.3cm, angle=90]{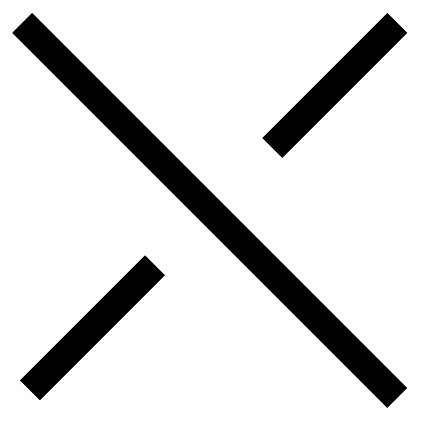}, we obtain
the $A$--smoothing as $\asplice$ and the $B$--smoothing as $\asymp$.
The state $s$ of $D$ is a choice of smoothing at every crossing,
resulting in a disjoint union of circles on $S^2$.  Let $|s|$ denote
the number of circles in $s$.  Let $s_A$ be the all--$A$ state, for
which every crossing of $D$ has an $A$--smoothing.  Similarly, $s_B$
is the all--$B$ state.  We will construct the Turaev surface
from the dual states $s_A$ and $s_B$.

\begin{figure}
\centerline{\includegraphics[height=85pt]{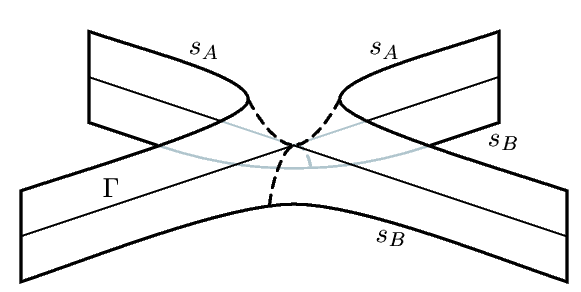}\ \ 
\includegraphics[height=75pt]{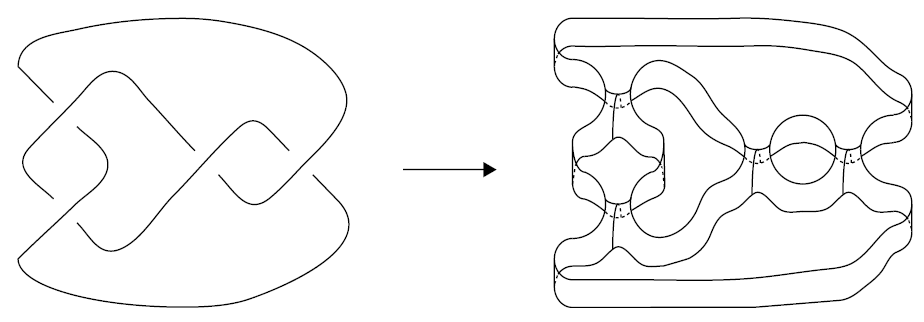}}
\caption{Cobordism between $s_A$ and $s_B$ (figures from \cite{cromwell} and \cite{Abe1})}
\label{figure:turaev-surface}
\end{figure}

At every crossing of $D$, we put a saddle surface which bounds the
$A$--smoothing on the top and the $B$--smoothing on the bottom as
shown in Figure \ref{figure:turaev-surface}.  In this way, we get a
cobordism between $s_A$ and $s_B$, with the link projection $\Gamma$ at the
level of the saddles.  The {\em Turaev surface} $F(D)$ is obtained by
attaching $|s_A|+|s_B|$ discs to all boundary circles.  See
Figure~\ref{figure:turaev-surface}, and \cite{hajij} for an animation
of the Turaev surface for the Borromean link.

The {\em Turaev genus} of $D$ is the genus of $F(D)$, given by
$$g_T(D) = g(F(D)) = (c(D) + 2 - |s_A| - |s_B|)/2.$$
The Turaev genus $g_T(L)$ of any non-split link $L$ is the minimum of $g_T(D)$ among all diagrams $D$ of $L$.
By \cite{turaev, DFKLS1}, $L$ is alternating if and only if $g_T(L)=0$, and if $D$ is an alternating diagram then $F(D) = S^2$.
In general, for any link diagram $D$, it follows that (see \cite{DFKLS1}):
\begin{enumerate}
\item $F(D)$ is a Heegaard surface of $S^3$; i.e., an unknotted closed orientable surface in $S^3$.
\item $D$ is alternating on $F(D)$, and the faces of $D$ can be
  checkerboard colored on $F(D)$, with discs for $s_A$ and $s_B$
  colored white and black, respectively.
\item $F(D)$ has a Morse decomposition, with $D$ and crossing saddles at height zero, and the $|s_A|$ and
  $|s_B|$ discs as maxima and minima, respectively.
\end{enumerate}

Conversely, in \cite{ADK} conditions were given for a Heegaard surface
with cellularly embedded alternating diagram on it to be a Turaev
surface.

\section{The Turaev surface and the Jones polynomial}\label{sec:Jones}

A diagram $D$ is $A$--adequate if at each crossing, the two arcs of
$s_A$ from that crossing are in different state circles.  In other
words, $|s_A|>|s|$ for any state $s$ with exactly one $B$--smoothing.
Similarly, we define a $B$--adequate diagram by reversing the roles of
$A$ and $B$ above.  If $D$ is both $A$--adequate and $B$--adequate it
is called {\em adequate}.  If $D$ is neither $A$--adequate nor
$B$--adequate it is called {\em inadequate}.  A link $L$ is {\em
  adequate} if it has an adequate diagram, and is {\em inadequate} if
all its diagrams are inadequate.  Any reduced alternating diagram is
adequate, hence every alternating link is adequate.  

Adequacy implies that $s_A$ and $s_B$ contribute the extreme terms
$\pm t^\alpha$ and $\pm t^\beta$ of the Jones polynomial $V_L(t)$,
which determine the span~$V_L(t) = |\alpha - \beta|$, which is a link
invariant.  Let $c(L)$ be the minimal crossing number among all
diagrams for $L$.  In \cite{turaev}, Turaev proved
$$ {\rm span}\, V_L(t) \leq c(L) - g_T(L) $$
with equality if $L$ is adequate.
If $D$ is a prime non-alternating diagram, then $g_T(D)>0$.
Thus, span~$V_L(t)= c(L)$ if and only if $L$ is alternating, from which Tait's Conjecture follows.

Therefore, for any adequate link $L$ with an adequate diagram $D$ (see \cite{Abe1}),
$$ g_T(L) = g_T(D) = \frac{1}{2}\left(c(D)-|s_A(D)|-|s_B(D)|\right)+1 = c(L) - {\rm span}\, V_L(t).$$
So for the connect sum $L\# L'$ of adequate links, $g_T(L\# L')=g_T(L)+g_T(L')$.

\ 

\noindent
\textbf{Turaev genus and knot homology.}
Khovanov homology and knot Floer homology categorify the Jones
polynomial and the Alexander polynomial, respectively.  The {\em
  width} of each bigraded knot homology, $w_{KH}(K)$ and $w_{HF}(K)$,
is the number of diagonals with non-zero homology.  The Turaev genus
bounds the width of both knot homologies \cite{dkh,lowrance1}:
\begin{equation}\label{KHineq} 
w_{KH}(K)-2 \leq g_T(K) \qquad \text{and} \qquad w_{HF}(K)-1 \leq g_T(K).
\end{equation}
For adequate knots, $w_{KH}(K)-2 = g_T(K)$ \cite{Abe1}.
These inequalities have been used to obtain families of knots with unbounded
Turaev genus (see \cite{CKsurvey}).

\ 

\noindent
\textbf{Ribbon graph invariants.}
Like the Seifert surface, the Turaev surface provides much more
information than its genus.  An oriented {\em ribbon graph} is a
graph with an embedding in an oriented surface, such that its faces
are discs.  Turaev's construction determines an oriented ribbon graph
$G_A$ on $F(D)$: We add an edge for every crossing in $s_A$, and
collapse each state circle of $s_A$ to a vertex of $G_A$, preserving
the cyclic order of edges given by the checkerboard coloring (see
\cite{CKsurvey}).

\pagebreak

If $L$ is alternating, then $V_L(t)=T_G(-t,-1/t)$, where $T_G(x,y)$ is
the Tutte polynomial \cite{Thistlethwaite}.  For any $L$, $V_L(t)$ is
a specialization of the Bollob\'as--Riordan--Tutte polynomial of $G_A$
\cite{DFKLS1}.  These ideas extend to virtual links and non-orientable
ribbon graphs \cite{Chmutov}.  In \cite{MMbook}, a unified description is
given for all these knot and ribbon graph polynomial invariants.

\begin{figure}
  \centerline{\includegraphics[height=100pt]{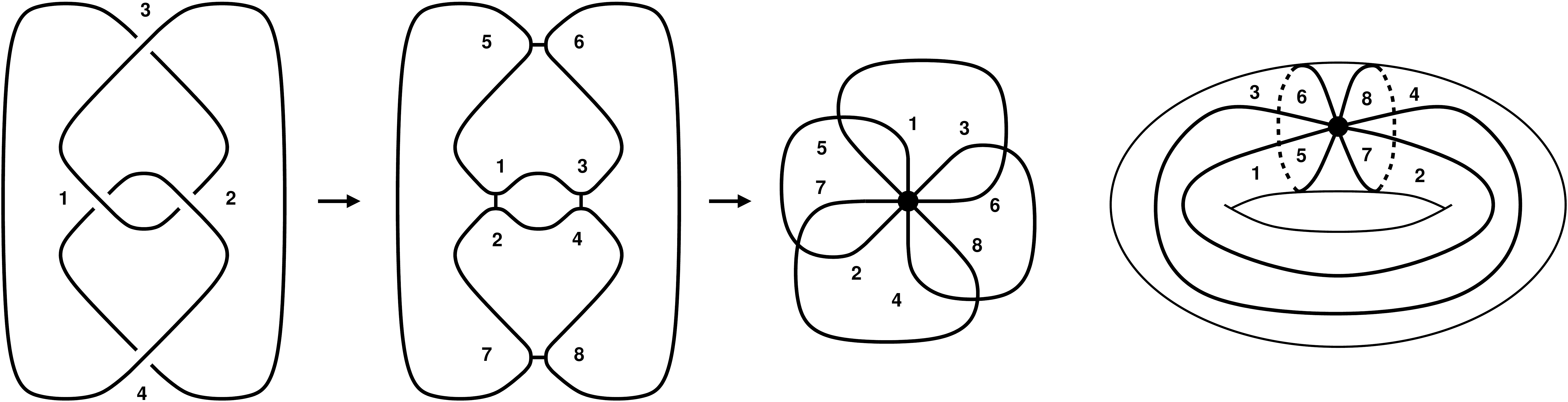}}
\caption{Ribbon graph $G_A$ for an inadequate diagram of the trefoil}
\label{figure:GA}
\end{figure}

\section{Turaev genus one links}\label{sec:g1}
The Seifert genus is directly computable for alternating and positive
links, and has been related to many classical invariants.  Moreover,
knot Floer homology detects the Seifert genus of knots. In contrast,
for most non-adequate links, computing the Turaev genus is an open
problem.

The Turaev genus of a link can be computed when the upper bounds in
the inequalities (\ref{KHineq}) or those in
\cite{Lowrance-knot_encyclopedia} match the Turaev genus of a
particular diagram, which gives a lower bound.  So it is useful to
know which diagrams realize a given Turaev genus.  Link diagrams with
Turaev genus one and two were classified in \cite{armond2017turaev, kim2018link}.

\begin{wrapfigure}{R}{0.3\textwidth}
\centering
\includegraphics[width=0.3\textwidth]{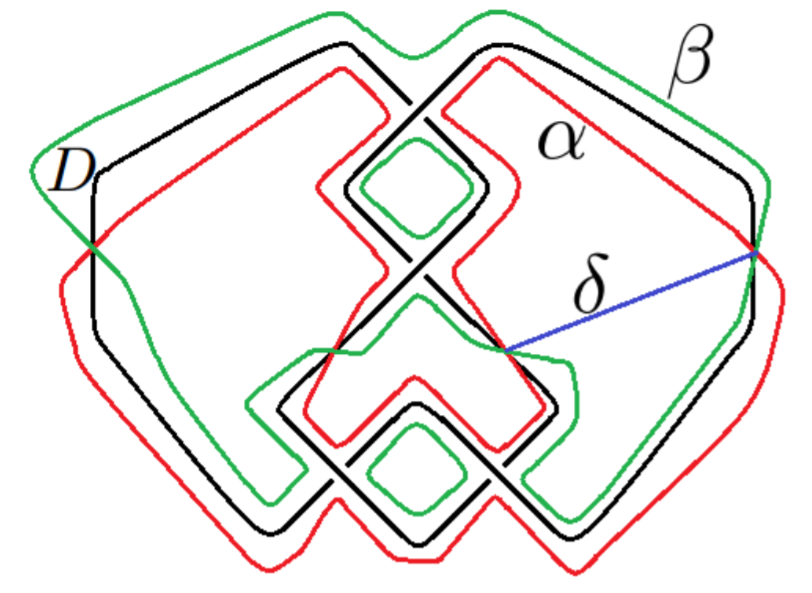}
\end{wrapfigure}

This classification uses the decomposition of any prime, connected
link diagram $D\subset S^2$ into alternating tangles.  An edge in $D$
is {\em non-alternating} when it joins two overpasses or two
underpasses.  If $D$ is non-alternating, we can isotope the state
circles in $s_A$ and $s_B$ to intersect exactly at the midpoints of
all non-alternating edges of $D$.  In the figure to the right from
\cite{kim2018link}, $\alpha\in s_A,\, \beta\in s_B$.  The arc $\delta$ joining
the points in $\alpha\cap\beta$ is called a {\em cutting arc} of $D$.

A cutting arc is the intersection of $S^2$ with a compressing disc of
the Turaev surface $F(D)$, which intersects $D$ at the endpoints of
$\delta$.  The boundary $\gamma$ of this compressing disc is called a
{\em cutting loop}.  Every cutting arc of $D$ has a corresponding
cutting loop on $F(D)$, and surgery of $D$ along a cutting arc
corresponds to surgery of $F(D)$ along a compressing disc, as shown in
the following figure from \cite{kim2018link}.

\pagebreak

\begin{figure}[h]
  \centerline{\includegraphics[height=110pt]{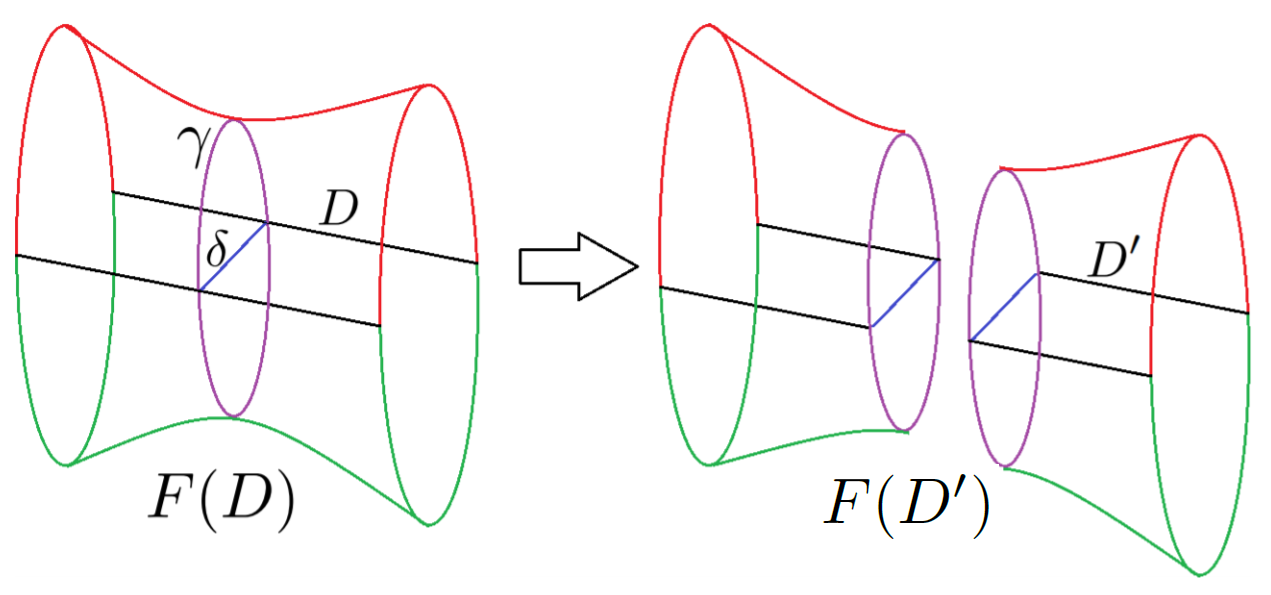}} 
\end{figure}

If $D'$ is obtained by surgery from $D$, the surgered surface is its
Turaev surface $F(D')$ with genus $g_T(D')=g_T(D)-1$.  So if
$g_T(D)=1$, then $\gamma$ is a meridian of the torus $F(D)$, and
surgery along all cutting arcs of $D$ cuts the diagram into
alternating $2$-tangles \cite{kim2018link}.  Hence, if $g_T(D)=1$,
then $D$ is a cycle of alternating $2$-tangles:

\begin{figure}[h]
  \centerline{\includegraphics[height=60pt]{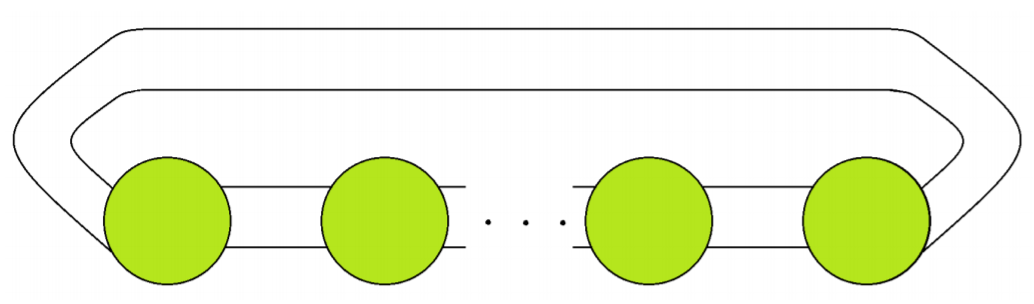}}
\end{figure}

\begin{wrapfigure}{R}{0.3\textwidth}
\centering
\includegraphics[width=0.3\textwidth]{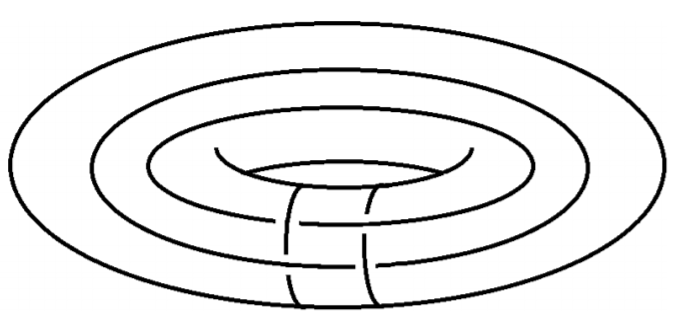}
\end{wrapfigure}

This also implies that for any alternating diagram $D$ on its Turaev
surface $F(D)$, if $g_T(D)\geq 1$ there is an essential simple loop
$\gamma$ on $F$ which intersects $D$ twice and bounds a disc in a
handlebody bounded by $F$.  Thus, the link on the surface in
Example~1.3.1 of \cite{Lowrance-knot_encyclopedia} cannot come from
Turaev's construction.  However, this condition is not sufficient; for
example, the diagram at right satisfies the condition, but cannot be a
Turaev surface because any planar diagram $D$ for this link has more
than four crossings, which would remain as crossings on $F(D)$.

Hayashi \cite{hayashi} and Ozawa \cite{ozawa2006non} considered more
general ways to quantify the complexity of the pair $(F,D)$, which has
prompted recent interest in {\em representativity} of knots (see,
e.g., \cite{arandarepresentativity, blair2017height,
  howie2017geometry, kindred2017alternating, ozawa2012bridge,
  ozawa2012representativity, pardon2011distortion}).

\section{Open problems}

Below, we consider open problems in two broad categories:

\begin{question}
How do you determine the Turaev genus of a knot or link?
\end{question}

Does the Turaev genus always equal the dealternating number of a link?
This is true in many cases, and no lower bounds are known to distinguish these invariants (see \cite{Lowrance-knot_encyclopedia}).

The lower bounds (\ref{KHineq}) vanish for quasi-alternating links.
For any $g>1$, does there exist a quasi-alternating link with Turaev genus $g$?

The Turaev genus is additive under connect sum for adequate knots, and
invariant under mutation if the diagram is adequate \cite{Abe1}.  In
general, for any $K$ and $K'$, is $g_T(K\#K')=g_T(K)+g_T(K')$?  If $K$
and $K'$ are mutant knots, is $g_T(K)=g_T(K')$?  The latter question
is open even for adequate knots; if $D$ is a non-adequate diagram of
an adequate knot $K$, then for a mutant $D'$ of $D$, it might be
possible that $g_T(K)<g_T(D)=g_T(D')=g_T(K')$.

If $K$ is a positive knot with Seifert genus $g(K)$, then $g_T(K)\leq g(K)$.
Is this inequality strict; i.e., is $g_T(K)<g(K)$ for a positive knot?
It is known to be strict for $g(K)=1,2$ \cite{jong2013positive, stoimenow2001knots} and for adequate positive knots \cite{lee2001adequate}.

In general, how do you compute the Turaev genus, which is a link invariant, without using link diagrams?
Is it determined by some other link invariants?

\begin{question}
How do you characterize the Turaev surface?
\end{question}

From the construction in Section~\ref{sec:def}, it is hard to tell
whether a given pair $(F,D)$ is a Turaev surface.  The existence of a
cutting loop implies that the alternating diagram $D$ on $F$ must have
minimal complexity; i.e., there exists an essential simple loop on $F$
which intersects $D$ twice.  But this condition is not sufficient.
What are the sufficient conditions for a given pair $(F,D)$ to be a
Turaev surface?  

Alternating, almost alternating and toroidally alternating knots have
been characterized topologically using a pair of spanning surfaces in
the knot complement \cite{greene2017alternating,
  howie2017characterisation, ito2018characterization,
  kim2016topological}.  Turaev genus one knots are toroidally
alternating, and they contain almost-alternating knots, but they have
not been characterized topologically as a separate class of knots.
What is a topological characterization of Turaev genus one knots, or
generally, of knots with any given Turaev genus?

Any non-split, prime, alternating link in $S^3$ is hyperbolic, unless
it is a closed 2--braid \cite{menasco}.  This result was recently
generalized to links in a thickened surface $F\times I$.  If the link
$L$ in $F\times I$ admits a diagram on $F$ which is alternating,
cellularly embedded, and satisfies an appropriate generalization of
``prime,'' then $(F\times I)-L$ is hyperbolic
\cite{adams2018hyperbolicity, ckp:bal_vol, howie2017geometry}.  Now,
for a given Turaev surface $F(D)$, let $L$ be a link in $F(D)\times I$
which projects to the alternating diagram on $F$.  It follows that
typically the complement $(F(D)\times I)-L$ is hyperbolic, assuming
there are no essential annuli. If $g_T(D)=1$ then $(F(D)\times I)-L$
has finite hyperbolic volume.  If $g_T(D)>1$ then there is a
well-defined finite volume if the two boundaries are totally geodesic.
How do the geometric invariants of $(F(D)\times I)-L$ depend on the
original diagram $D$ in $S^2$?

\bibliographystyle{plain}
\bibliography{Turaev_surfaces_references}

\end{document}